\documentclass[preprint,3p,12pt,times]{elsarticle} 
\usepackage{amssymb}
\usepackage{amsthm}
\usepackage{amsmath}

\newtheorem{thm}{Theorem}

\newtheorem{cor}[thm]{Corollary}
\newtheorem{lem}[thm]{Lemma}
\newtheorem*{open}{Open Problem}

\newtheorem{knownthm}{Theorem}

\newtheorem{knownlem}[knownthm]{Lemma}

\theoremstyle{definition}

\newcommand{\C}{\mathbb{C}}
\newcommand{\D}{\mathbb{D}}

\newcommand{\A}{\mathcal{A}}
\renewcommand{\S}{\mathcal{S}}
\newcommand{\K}{\mathcal{K}}
\newcommand{\pf}{\noindent{\bf\textit{Proof. \,}}}
\newcommand{\epf}{\phantom{} \hfill$ \Box$}
\newcommand{\no}{\noindent}
\newcommand{\dstyle}{\displaystyle}
\renewcommand{\Re}{\textup{Re}\,}

\newcommand{\Tan}{\rm Tan}

\begin{document}


\begin{frontmatter}

\title{On strongly starlike and convex functions of order $\alpha$ and type $\beta$} 

\author{IKKEI HOTTA\fnref{label2}}
\ead{ikkeihotta@ims.is.tohoku.ac.jp} 
\address{Division of Mathematics, Graduate School of Information Sciences, Tohoku University, 6-3-09 Aramaki-Aza-Aoba, Aoba-ku, Sendai, Miyagi 980-8579, Japan \\Telephone number ; +81-22-795-4636} 
\fntext[label2]{Corresponding author.}

\author{MAMORU NUNOKAWA} 
\ead{mamoru\_nuno@doctor.nifty.jp} 
\address{Gunma University, Hoshikuki 798-8, Chuou-Ward, Chiba 260-0808, Japan}

\begin{abstract}
In this note we investigate the inclusion relationship between the class of strongly starlike functions of order $\alpha$ and type $\beta$, $\alpha \in (0,1]$ and $\beta \in [0,1)$, which satisfy
$$
\left|
 \arg
 \left\{
  \frac{zf'(z)}{f(z)} - \beta
 \right\}
\right|
<
\frac{\pi}{2}\alpha
$$
and the class of strongly convex functions of order $\alpha$ and type $\beta$ which satisfy
$$
\left|
 \arg
  \left\{
   1 + \frac{zf''(z)}{f'(z)} -\beta
  \right\}
\right| 
<
\frac{\pi}{2}\alpha
$$
in the unit disk, where $f$ is an analytic function defined on the unit disk and satisfies $f(0)=f'(0)-1 = 1$.
Some applications of our main result are also presented which contains various classical results for the typical subclasses of starlike and convex functions.
\end{abstract}


\end{frontmatter}


\section{Introduction}
Let $\A$ denote the set of functions $f(z) = z + \sum_{n=2}^{\infty}a_{n}z^{n}$ which are analytic in the unit disk $\D = \{z \in \C : |z|<1\}$.

Let $\alpha$ be a real number with $\alpha \in (0,1]$.
A function $f \in \A$ is called \textit{strongly starlike of order $\alpha$} if it satisfies
$$
\left|
 \arg 
  \left\{
   \frac{zf'(z)}{f(z)}
  \right\}
\right|
< \frac{\pi}{2}\alpha
$$
for all $z \in \D$ and \textit{strongly convex of order $\alpha$} if
$$
\left|
 \arg 
  \left\{
   1 + \frac{zf''(z)}{f'(z)}
  \right\}
\right|
< \frac{\pi}{2}\alpha
$$
for all $z \in \D$.
Let us denote by $\S^{*}(\alpha)$ the class of functions strongly starlike of order $\alpha$, and by $\K(\alpha)$ the class of functions strongly convex of order $\alpha$.
The class $\S^{*}(\alpha)$ was introduced first by Stankiewicz \cite{Stankiewicz:1966} and by Brannan and Kirwan \cite{BrannanKirwan:1969}, independently.
It is clear from the definitions that $\S^{*}(\alpha_{1}) \subset \S^{*}(\alpha_{2})$ and $\K(\alpha_{1}) \subset \K(\alpha_{2})$ for $0 < \alpha_{1} < \alpha_{2} \leq 1$.
The case when $\alpha = 1$, i.e., $\S^{*}(1)$ and $\K(1)$ correspond to well known classes of starlike and convex functions respectively, and therefore all the functions which belong to $\S^{*}(\alpha)$ or $\K(\alpha)$ are univalent in $\D$.
We denote by $\S^{*}$ and $\K$ the classes of starlike and convex functions.
For the general reference of classes of starlike and convex functions, see, for instance \cite{Duren:1983}.

Mocanu \cite{Mocanu:1989} obtained the following result (see also \cite{Nunokawa:1993a}).
Here, set
\begin{equation}\label{rho}
\rho(\alpha) = \Tan^{-1}
\frac
 {
  \dstyle 
  \left(\frac{\alpha}{1 - \alpha}\right)
  \left(\frac{1-\alpha}{1+\alpha}\right)^{\frac{1+ \alpha}{2}} 
  \sin\left[\frac{\pi}{2}(1-\alpha)\right]
 }
 {
  \dstyle
  1 + 
  \left(\frac{\alpha}{1 - \alpha}\right)
  \left(\frac{1-\alpha}{1+\alpha}\right)^{\frac{1+ \alpha}{2}} 
  \cos\left[\frac{\pi}{2}(1-\alpha)\right]
 }
\end{equation}
and
$$
\gamma(\alpha) = \alpha + \frac{2}{\pi}\rho(\alpha).
$$

\begin{knownthm}\label{mocanu}
$\K(\gamma(\alpha)) \subset \S^{*}(\alpha)$ for each $\alpha \in (0,1]$.
\end{knownthm}

We remark that the function $\gamma(\alpha)$ is continuous and strictly increases from 0 to 1 when $\alpha$ moves from 0 to 1.
Further investigations for the above theorem can be found in \cite{KanasSugawa:2007}.

Now we shall introduce the class of functions $\S^{*}(\alpha,\beta)$ and $\K(\alpha,\beta),\,\alpha \in (0,1]$ and $\beta \in [0,1)$, whose members satisfy the conditions
$$
\left|
 \arg
 \left\{
  \frac{zf'(z)}{f(z)} - \beta
 \right\}
\right|
<
\frac{\pi}{2}\alpha
$$
and
$$
\left|
 \arg
  \left\{
   1 + \frac{zf''(z)}{f'(z)} -\beta
  \right\}
\right| 
<
\frac{\pi}{2}\alpha
$$
for all $z \in \D$, respectively.
We call a function $f \in \S^{*}(\alpha,\beta)$ \textit{strongly starlike of order $\alpha$ and type $\beta$}. 
In the same way, a function $f \in \K(\alpha,\beta)$ is \textit{strongly convex of order $\alpha$ and type $\beta$}.
It is obvious that $\S^{*}(\alpha,0) = \S^{*}(\alpha)$ and $\K(\alpha,0) = \K(\alpha)$.
Also the following relations are true from the definitions; 
\begin{equation}\label{relations}
\begin{array}{rlll}
\textup{i)}&\S^{*}(\alpha_{1},\beta) \subset \S^{*}(\alpha_{2},\beta),\hspace{200pt}\\
\textup{ii)}&\K(\alpha_{1},\beta) \subset \K(\alpha_{2},\beta),\\
\textup{iii)}&\S^{*}(\alpha,\beta_{1}) \supset \S^{*}(\alpha,\beta_{2}),\\
\textup{iv)}&\K(\alpha,\beta_{1}) \supset \K(\alpha,\beta_{2}),
\end{array}
\end{equation}
for $0 < \alpha_{1} < \alpha_{2} \leq 1$ and $0 \leq \beta_{1} < \beta_{2} < 1$.
That is why all functions belong to $\S^{*}(\alpha,\beta)$ or $\K(\alpha,\beta)$ are univalent on $\D$.

A sufficient condition for which $f \in \A$ lies in $\S^{*}(\alpha, \beta)$ was proved by the second author et al.~\cite{Nunokawa.etal:1997a}.
The authors also proposed in \cite{Nunokawa.etal:1997a} the open problem about a inclusion relationship between $\K(\alpha,\beta)$ and $\S^{*}(\alpha,\beta)$.
However, it seems that no results concerning this question have been known.

Our main result is the following;

\begin{thm}\label{main}
$\K(\gamma(\alpha),\beta) \subset \S^{*}(\alpha,\beta)$ for each $\alpha \in (0,1]$ and $\beta \in [0,1)$.
\end{thm}

The above theorem includes Theorem \ref{mocanu} as the case when $\beta=0$.

We should notice the reader that this estimation is not sharp for each $\alpha \in (0,1]$ and $\beta \in [0,1)$ (see also \cite{KanasSugawa:2007}).
We will discuss about this problem in section 2 with the proof of Theorem \ref{main}.
Our main theorem yields several applications which will be shown in the last section.


\section{Proof of Theorem \ref{main}}\label{proof}

Our proof relies on the following lemma which was obtained by the second author \cite{Nunokawa:1992a,Nunokawa:1993a};

\begin{knownlem}\label{Nunokawa}
Let $p(z)$ be analytic and satisfies $p(0)=1,\,p(z) \neq 0$ in $\D$.
Let us assume that there exists a point $z_{0} \in \D$ such that $|\arg p(z)| < \pi\alpha/2$
for $|z| < |z_{0}|$ and $|\arg p(z_{0})| = \pi\alpha/2$
where $\alpha >0$.
Then we have
$$
\frac{z_{0}p'(z_{0})}{p(z_{0})} = i \alpha k
$$
where
$
k \geq \frac12\left(a + \frac1a\right)
$
when $\arg p(z_{0}) = \pi\alpha/2$ and
$
k \leq - \frac12\left(a + \frac1a\right)
$
when $\arg p(z_{0}) = - \pi\alpha/2$, where $p(z_{0})^{1/\alpha} = \pm i a$ and $a>0$.
\end{knownlem}

The next result will be used later;
\begin{lem}\label{lem2}
$\Tan^{-1}\alpha \geq \rho(\alpha)$ for all $\alpha \in (0,1]$, where $\rho$ is defined by \eqref{rho}.
\end{lem}

\pf
Put
$
\phi(\alpha) = 
\left(1/(1 - \alpha)\right)
\left((1-\alpha)/(1+\alpha)\right)^{\frac{1+ \alpha}{2}}.
$
It is enough to prove that
$$
\alpha \geq
\frac
 {
  \dstyle 
  \alpha \phi(\alpha)
  \sin\left[\pi(1-\alpha)/2\right]
 }
 {
  \dstyle
  1 + 
  \alpha \phi(\alpha)
  \cos\left[\pi(1-\alpha)/2\right]
 }
$$
for all $\alpha \in (0,1]$.
Since $\phi(\alpha) <1$ because of $\phi(0)=1$ and $\phi'(\alpha) <0$, we obtain $\alpha > \phi(\alpha)$ and therefore
$$
\frac
 {
  \dstyle 
  \alpha
  \sin\left[\pi(1-\alpha)/2\right]
 }
 {
  \dstyle
  1 + 
  \alpha
  \cos\left[\pi(1-\alpha)/2\right]
 }
>
\frac
 {
  \dstyle 
  \alpha\phi(\alpha)
  \sin\left[\pi(1-\alpha)/2\right]
 }
 {
  \dstyle
  1 + 
  \alpha\phi(\alpha)
  \cos\left[\pi(1-\alpha)/2\right]
 }.
$$
It remains to show that 
$$
\alpha \geq
\frac
 {
  \dstyle 
  \alpha
  \sin\left[\pi(1-\alpha)/2\right]
 }
 {
  \dstyle
  1 + 
  \alpha
  \cos\left[\pi(1-\alpha)/2\right]
 }
$$
for all $\alpha \in (0,1]$ and this is clear.
\epf

\

\no
\textbf{\textit{Proof of Theorem \ref{main}.\,}}
Let us suppose that $f$ satisfies the assumption of the theorem and let
$$
p(z) = 
\frac{1}{1-\beta}
\left(
 \frac{zf'(z)}{f(z)} - \beta
\right).
$$
Then $p(0) =1$, and calculations show that
\begin{equation}\label{proof02}
1 + 
\frac{zf''(z)}{f'(z)} - \beta
=
(1-\beta) p(z)
\left\{  
 1 + 
 \frac
 {\frac{zp'(z)}{p(z)}}
 { (1-\beta)p(z) + \beta}
\right\}.
\end{equation}
We note that $p(z) \neq 0$ holds for all $ z \in \D$ since $1 + zf''(z)/f'(z) - \beta \neq \infty$ on $\D$ from our assumption.

Now we derive a contradiction by using Lemma \ref{Nunokawa}.
If there exists a point $z_{0}$ such that $|\arg p(z)| < \pi\alpha/2$
for $|z| < |z_{0}|$ and $|\arg p(z_{0})| = \pi\alpha/2$, where $\alpha \in (0,1]$, then by Lemma \ref{Nunokawa}, $p$ must satisfy
$
z_{0}p'(z_{0})/p(z_{0}) = i \alpha k
$
where
$
k \geq \frac12\left(a + \frac1a\right)
$
when $\arg p(z_{0}) = \pi\alpha/2$ and
$
k \leq - \frac12\left(a + \frac1a\right)
$
when $\arg p(z_{0}) = - \pi\alpha/2$, where $p(z_{0})^{1/\alpha} = \pm i a$ and $a>0$.

At first we suppose that $\arg p(z_{0}) = \pi\alpha/2$.
Then from \eqref{proof02} we have
\begin{eqnarray*}
\arg
\left\{
 1 + \frac{z_{0}f''(z_{0})}{f'(z_{0})} - \beta
\right\}
&=&
\arg
\left[
 (1-\beta) p(z_{0})
\left\{  
 1 + 
 \frac
 {\frac{z_{0}p'(z_{0})}{p(z_{0})}}
 { (1-\beta)p(z_{0}) + \beta}
\right\}
\right]\\[10pt]
&=&
\frac{\pi}{2}\alpha + 
\arg
\left\{  
 1 + 
 \frac
 {i \alpha k}
 { (1-\beta)p(z_{0}) + \beta}
\right\}.
\end{eqnarray*}

\no
We shall estimate the second term of the second line of above.
Geometric observations show that the point $1 + [i \alpha k / \{(1 - \beta)p(z_{0}) +\beta\}]$ lies on the subarc $C$ of the circle which passes through $1,\,1 + i \alpha k$ and $1 + [i \alpha k / p(z_{0})]$, where $C$ connects $1 + i \alpha k$ and $1 + [i \alpha k /p(z_{0})]$ and does not pass through 1.
Further, we can find out that the value $\{\arg z : z \in C\}$ attains its minimum at the end points of $C$.
Therefore we have
\begin{equation}\label{inequality1}
\arg
\left\{  
 1 + 
 \frac
 {i \alpha k}
 { (1-\beta)p(z_{0}) + \beta}
\right\}
\geq
\min 
\left\{
 \arg
 \left\{
  1 + i\alpha k
 \right\},\,
 \arg
 \left\{
  1 + \frac{i\alpha k}{p(z_{0})}
 \right\}
\right\}.
\end{equation}

Here, the first value in the above minimum can be evaluated by $\arg \{1 + i \alpha k\} \geq \Tan^{-1} \alpha$ since $k \geq 1$.
For the second value, we note that $a^{1-\alpha} + a^{-1-\alpha}$ takes its minimum value at $a = \sqrt{(1+\alpha)/(1-\alpha)}$.
Therefore
\begin{eqnarray*}
\arg
\left\{
 1 + \frac{i \alpha k}{p(z_{0})} 
\right\}
&=&
\arg
\left\{
 1 + e^{\frac{\pi}{2}(1-\alpha)i}
 \cdot
 \frac{\alpha}{2}
 [a^{1-\alpha} + a^{-1-\alpha}]
\right\}\\[8pt]
&\geq&
\arg
\left\{
 1 + e^{\frac{\pi}{2}(1-\alpha)i}
 \cdot
 \frac{\alpha}{2}
 \left[
  \left(
   \frac{1 + \alpha}{1 - \alpha}
  \right)^{\frac{1-\alpha}{2}}
  +
  \left(
   \frac{1 + \alpha}{1 - \alpha}
  \right)^{\frac{-1-\alpha}{2}}
 \right]
\right\}\\[8pt]
&=&
\rho(\alpha).
\end{eqnarray*}

By Lemma \ref{lem2} we conclude that
$$
\arg
\left\{
 1 + \frac{z_{0}f''(z_{0})}{f'(z_{0})} - \beta
\right\}
\geq
\frac{\pi}{2}\alpha + 
\min
 \left\{
  \Tan^{-1}\alpha,\,\rho(\alpha)
 \right\} 
= 
\frac{\pi}{2}\gamma(\alpha)
$$
and this contradicts our assumption.

In the same fashion, if $\arg p(z_{0}) = - \pi \alpha /2$ then a similar argument shows that
$$
\arg
\left\{
 1 + \frac{z_{0}f''(z_{0})}{f'(z_{0})} - \beta
\right\}
\leq
-\frac{\pi}{2}\alpha + 
\max
 \left\{
  \Tan^{-1}(-\alpha),\,-\rho(\alpha)
 \right\} 
= 
-\frac{\pi}{2}\gamma(\alpha).
$$
This also contradicts our assumption and our proof is completed.
\epf

We remark that we expect this theorem to be room for improvement in our method because the inequality \eqref{inequality1} is a rough estimation except the case when $\beta = 0$, whereas it seems to be not easy to give a precise estimation for the left hand side of \eqref{inequality1}.


\section{Applications}\label{appli}

We would like to give a further discussion to the relationship between $\S^{*}(\alpha,\beta)$ and $\K(\alpha,\beta)$ by using Theorem \ref{main}.

\subsection{}
It is well known that a convex function is a starlike function, that is, $\K \subset \S^{*}$.
Furthermore, Mocanu \cite{Mocanu:1986a} showed that $\K(\alpha) \subset \S^{*}(\alpha)$ for all $\alpha \in (0,1]$. 
Now we give the next result which includes these properties as special cases;

\begin{cor}\label{cor1}
$\K(\alpha,\beta) \subset \S^{*}(\alpha,\beta)$ for each $\alpha \in (0,1]$ and $\beta \in [0,1)$.
\end{cor}
\pf
Since $\alpha \leq \gamma(\alpha)$ for all $\alpha \in (0,1]$, $\K(\alpha,\beta) \subset \K(\gamma(\alpha),\beta) \subset \S^{*}(\alpha,\beta)$ by ii) in \eqref{relations} and Theorem \ref{main} which is our desired inclusion.
\epf

Corollary \ref{cor1} yields the following property;
\begin{cor}
If $zf'(z) \in \S^{*}(\alpha,\beta)$, then $f \in \S^{*}(\alpha,\beta)$.
\end{cor}
\pf
It is obvious that $g \in \K(\alpha,\beta)$ if and only if $zg'(z) \in \S^{*}(\alpha,\beta)$.
Thus if $zg'(z) \in \S^{*}(\alpha,\beta)$ then $g \in \K(\alpha,\beta) \subset \S^{*}(\alpha,\beta)$ from Corollary \ref{cor1}.
Hence our assertion follows if we put $f(z)=zg'(z)$.
\epf

This corollary is equivalent to the following;
\textit{$\S^{*}(\alpha,\beta)$ is preserved by the Alexander transformation,
where the Alexander transformation \cite{Alexander:1915} is the integral transformation defined by
$$
f(z) \mapsto \int_{0}^{z}\frac{f(u)}{u}du
$$
for $f \in \A$.
}

\subsection{}
If $\alpha = 1$, then the class $\S^{*}(1,\beta)$ and $\K(1,\beta)$ is called \textit{starlike of order $\beta$} and \textit{convex of order $\beta$}, respectively.
It is easy to see that $f \in \S^{*}(1,\beta)$ satisfies
$$
\Re \left\{ \frac{zf'(z)}{f(z)} \right\} > \beta
$$
and $f \in \K(1,\beta)$ satisfies
$$
\Re \left\{ 1 + \frac{zf''(z)}{f'(z)} \right\} > \beta.
$$
Marx \cite{Marx:1933} and Strohh{\"a}cker \cite{Strohhacker:1933} showed that $\K(1,0) \subset \S^{*}(1,1/2)$.
Jack \cite{Jack:1971} proposed the more general problem; What is the largest number $\beta_{0}$ which satisfies $\K(1,\beta) \subset \S^{*}(1,\beta_{0})$?
Later MacGregor \cite{MacGregor:1975} and Wilken and Feng \cite{WilkenFeng:1980} answered the problem to give the exact value of $\beta_{0}$;

\begin{knownthm}\label{MacGregoretal}
$\K(1,\beta) \subset \S^{*}(1,\delta(\beta))$ for all $\beta \in [0,1)$, where
$$
\delta(\beta) =
\left\{
 \begin{array}{ccc}
 \dstyle\frac{1-2\beta}{2^{2-2\beta}(1-2^{2\beta-1})}, & 
 \textit{if} & 
 \dstyle \beta \neq \frac12, \\[10pt]
 \dstyle \frac{1}{2 \log 2}, &
 \textit{if} &
 \dstyle \beta = \frac12.
 \end{array}
\right.
$$
This estimation is sharp for each $\beta \in [0,1)$.
\end{knownthm}

Setting $\beta = 0$, we have the result of Marx and Strohh{\"a}cker.
We can obtain a similar estimation to above that ``$\K(\gamma(\alpha),\delta(\beta)) \subset \S^{*}(\alpha,\beta)$ for all $\alpha \in (0,1]$ and $\beta \in [0,1)$'' by Theorem \ref{main} since $\beta < \delta(\beta)$ for all $\beta \in [0,1)$. 
However, the following problem is still open;
\begin{open}
$\K(\gamma(\alpha),\beta) \subset \S^{*}(\alpha,\delta(\beta))$ for each $\alpha \in (0,1]$ and $\beta \in [0,1)$.
\end{open}
This problem implies Theorem \ref{main} because $\S^{*}(\alpha,\delta(\beta)) \subset \S^{*}(\alpha,\beta)$ for all $\alpha \in (0,1]$ and $\beta \in [0,1)$, and Theorem \ref{MacGregoretal} as the case when $\alpha = 1$.

\bibliographystyle{amsplain}
\bibliography{bibdata}

\end{document}